\documentclass[english,11pt,a4paper,leqno]{amsart}
\usepackage{amsmath,amstext, amsthm, amssymb, mathtools}
\usepackage{latexsym, amsthm, amsfonts, amsmath,amssymb,graphicx,xcolor,dsfont}
\usepackage[latin1]{inputenc} 
\usepackage{enumerate}
\usepackage{enumitem}
\usepackage{hyperref}
\usepackage[english]{babel} 

\hypersetup{
	colorlinks   = true, 
	urlcolor     = blue, 
	linkcolor    = blue, 
	citecolor   = red 
}

\newtheorem{thm}{Theorem}[section]

\newtheorem{lem}[thm]{Lemma}
\newtheorem{coro}[thm]{Corollary}

\newtheorem*{thm*}{Theorem}

\newtheorem{rem}[thm]{Remark}

\newtheorem*{nota*}{Notation}
\newtheorem*{defi*}{Definition}

\numberwithin{equation}{section}

\newcommand{\Q}{{\overline{\mathbb Q}}}
\newcommand{\C}{{\mathbb{C}}}
\newcommand{\Z}{{\mathbb{Z}}}
\newcommand{\bQ}{{\overline{\mathbb Q}}}
\newcommand{\N}{{\mathbb{N}}}

\newcommand{\I}{\mathcal I}

\newcommand{\z}{{\boldsymbol{z}}}
\newcommand{\X}{{\boldsymbol{X}}}

\newcommand{\bmu}{{\boldsymbol{\mu}}}

\newcommand{\bbeta}{{\boldsymbol{\beta}}}
\newcommand{\balpha}{{\boldsymbol{\alpha}}}

\newcommand{\M}{{\boldsymbol{\rm M}}}
\numberwithin{equation}{section} 

\title{Linear disjointness and values of Mahler $M$-functions}

\author{Boris Adamczewski}
\address{
    Univ Lyon, Universit\'e Claude Bernard Lyon 1\\
    CNRS UMR 5208, Institut Camille Jordan, F-69622 Villeurbanne Cedex, France,  and 
    Centre International de Rencontres Math\'ematiques (CIRM), CNRS UAR 822, 163 Av.\ de Luminy, 13009
    Marseille, France}
\email{boris.adamczewski@math.cnrs.fr}

\author{Colin Faverjon}
\address{
Universit\'e Picardie Jules Verne \\
UMR CNRS 7352, LAMFA \\
33, rue saint-Leu 80039 Amiens, France}
\email{colin.faverjon@math.cnrs.fr}
\date{}

\thanks{ 
Funded by the European Union (ERC Advanced Grant M\^{E}TIS, Grant Agreement No. 101264813). Views and opinions expressed are however those of the authors only and do not necessarily reflect those of the European Union or the European Research Council Executive Agency. Neither the European Union nor the granting authority can be held responsible for them.}
\begin{document}
\maketitle

\begin{abstract}
We prove that the algebra generated by the values of Mahler $M$-functions at
algebraic points admits a  decomposition as the tensor product,
over $\Q$, of a countable family of linearly disjoint algebras.
Consequently, the study of linear and algebraic relations among such
values reduces to the corresponding problem within each of these 
smaller algebras. Our proof relies on recent advances of the authors on Mahler's
method in several variables. As an application, we obtain a complete
description of the algebraic relations among the values of the Fredholm
series and their derivatives at nonzero algebraic points. This fully
resolves an old problem originating in Mahler's work.
\end{abstract}

\section{Introduction}\label{sec:intro}

A recurrent strategy throughout mathematics consists in decomposing a
complicated problem into simpler subproblems and then recovering the
solution of the original problem from the solutions of its components.
In representation theory, for instance, one seeks to decompose a
representation into irreducible representations. In dynamical systems,
one studies invariant subsets and, whenever possible, decomposes the
system into ergodic components. Spectral decompositions provide another
familiar manifestation of the same principle. In all these situations, the global structure can be recovered from
suitably chosen elementary pieces.

Questions of linear or algebraic independence are, however, intrinsically
global, which prevents such a strategy from working in general.  To illustrate
this elementary obstruction, let $x_1,\ldots,x_n$ be linearly
(resp.\ algebraically) independent over a field $\mathbb K$. Then every
proper subset of
\[
\{x_1,\ldots,x_n,x_1+\cdots+x_n\}
\]
is linearly (resp.\ algebraically) independent over $\mathbb K$, whereas
the whole set is not. Thus, even a complete understanding of the
relations within all proper subfamilies does not determine the relations
within the entire family.

In this paper, we exhibit what seems to be the first general
unconditional setting in transcendental number theory in which the
global obstruction described above can nevertheless be overcome. More precisely, we prove that the
$\Q$-algebra $\M$ generated by all values of Mahler $M$-functions at
algebraic points admits a  decomposition as the tensor product,
over $\Q$, of a countable family of linearly disjoint algebras. Consequently, the study of all linear and
algebraic relations among arbitrary $M$-values reduces to the
corresponding problem within each of these smaller algebras. 
This class of numbers has been studied since the pioneering work of
Mahler and has attracted renewed interest in recent years, notably
because of its close connections with finite automata and automatic
numbers, as well as with a conjecture of Furstenberg concerning
expansions of real numbers in multiplicatively independent bases; see
\cite{AF26,AF26_ad}.

Let $\mathbb K\subset\mathbb L$ be fields, and let $\mathbb E$ and
$\mathbb F$ be two $\mathbb K$-subalgebras of $\mathbb L$. They are said
to be \emph{linearly disjoint over $\mathbb K$} if the canonical map
\[
\begin{array}{rcl}
\mathbb E\otimes_{\mathbb K}\mathbb F&\longrightarrow&\mathbb L,\\
\xi\otimes\eta&\longmapsto&\xi\eta,
\end{array}
\]
is injective. Equivalently, every finite family of elements of $\mathbb E$
that is linearly independent over $\mathbb K$ remains linearly independent
over $\mathbb F$. Similarly, we say that $\mathbb E$ and $\mathbb F$ are
\emph{free over $\mathbb K$} if every finite family of elements of
$\mathbb E$ that is algebraically independent over $\mathbb K$ remains
algebraically independent over $\mathbb F$. When $\mathbb K$ is algebraically closed and has characteristic zero,
linear disjointness and freeness are equivalent.  
These notions extend naturally from two algebras to an arbitrary family.
Let $\mathbb K\subset\mathbb L$ be fields, and let
$(\mathbb E_i)_{i\in I}$ be a family of $\mathbb K$-subalgebras of
$\mathbb L$. We say that the algebras $\mathbb E_i$, $i\in I$, are
linearly disjoint (resp.\ free) over $\mathbb K$ if, for every
$i_0\in I$, the algebras
$$
\mathbb E_{i_0}
\quad\mbox{ and }\quad
\mathbb K[\mathbb E_i:i\in I\setminus\{i_0\}]
$$
are linearly disjoint (resp.\ free) over $\mathbb K$.  

Linear disjointness therefore provides precisely the decomposition
principle that is missing for arbitrary families of numbers: once the
relations within each individual algebra are known, all global relations
are obtained by combining in a trivial way the relations within each
subalgebra. 

It is a fundamental notion in field theory and Galois theory. Variants
of the same principle also play an important role in differential and
difference Galois theory, as well as in the model-theoretic study of
independence. By contrast, its explicit use in transcendental number
theory seems to be extremely scarce. The main examples we could find
occur in \cite{Cheng09,PSS20,BPSS}. They are conditional results relying
on Schanuel's conjecture, the Grothendieck period conjecture, or
Andr\'e's generalization of the latter, and were developed following a
remark of Lang \cite[pp.\ 30--31]{Lang_trans}. For instance, assuming
Schanuel's conjecture, the authors of \cite{Cheng09} proved that
\[
\Q[e,e^e,e^{e^e},e^{e^{e^e}},\ldots]
\quad\text{and}\quad
\Q[\pi,\log(\pi),\log(\log(\pi)),
\log(\log(\log(\pi))),\ldots]
\]
are linearly disjoint over $\Q$, where, at each stage, $\log$ denotes a
suitable determination of the logarithm.

\subsection{Linear disjointness for $M$-values}

We now state our main results more precisely. Throughout the paper,
we let $\Q$ denote the algebraic closure of $\mathbb Q$ in $\C$, and $\Q^\star=\Q\setminus\{0\}$. 

Let $q\geq 2$ be an integer. 
 An
\emph{$M_q$-function} is a power series $f\in\Q[[z]]$ for which there
exist polynomials $a_0,\ldots,a_m\in\bQ[z]$, not all zero, such that
\begin{equation}
\label{eq:mahler_linear}
a_0(z)f(z)+a_1(z)f(z^q)+\cdots+a_m(z)f(z^{q^m})=0.
\end{equation}
An \emph{$M$-function} is a power series that is an $M_q$-function for
some integer $q\geq 2$. We recall that an $M$-function is either rational
or transcendental over $\Q(z)$. In either case, it is analytic in a
neighborhood of the origin and meromorphic in the open unit disk. In the
transcendental case, the unit circle is a natural boundary. Consequently,
such a function can only be evaluated at points of the open unit disk $\mathcal D=\{\xi \in \C: \vert \xi \vert < 1\}$ 
that are not poles.

Our main object of interest is the $\Q$-algebra
$$
\M:=
\Q\bigl[
f(\alpha):
\alpha\in\Q^\star\cap\mathcal D,\ 
\text{$f$ is an $M$-function well-defined at $\alpha$}
\bigr].
$$

Let $\mathbb K\subset\Q$ be a field, let $q\geq 2$, and let
$\alpha\in\mathbb K$ satisfy $0<|\alpha|<1$. We set
$$
\M_{\mathbb K,q,\alpha}
:=
\bigl\{
f(\alpha):
f\in\mathbb K[[z]]
\text{ is an $M_q$-function well-defined at $\alpha$}
\bigr\}.
$$
This defines a family of $\mathbb K$-algebras.

Recall that two nonzero complex numbers $\theta$ and $\eta$ are
\emph{multiplicatively dependent} if
$\theta^u=\eta^v$ for some nonzero pair of integers $(u,v)$. Given
$(q,\alpha)$ and $(r,\beta)$ in $\mathbb Z_{\geq 2}\times\Q^\star$, we write
$(q,\alpha)\sim(r,\beta)$ if $q$ and $r$ are multiplicatively dependent
and $\alpha$ and $\beta$ are multiplicatively dependent. This defines an
equivalence relation.

Our main result reads as follows.

\begin{thm}
\label{thm:M-values}
Let $\mathbb K\subset\Q$ be a field.
 Then the
following assertions hold.
\begin{enumerate}[label=\rm (\alph*)]
\item\label{item:equal}
Let $(q,\alpha),(r,\beta)\in\mathbb Z_{\geq 2}\times(\mathbb K^\star\cap \mathcal D)$. Then
$$
\M_{\mathbb K,q,\alpha}
=
\M_{\mathbb K,r,\beta} \iff (q,\alpha)\sim(r,\beta).
$$

\item
Let $\I\subset \mathbb Z_{\geq 2}\times(\mathbb K^\star\cap \mathcal D)$ be such that
$(q,\alpha)\not\sim(r,\beta)$ for every pair of distinct elements
$(q,\alpha),(r,\beta)\in\I$. Then the family of $\mathbb K$-algebras
$\M_{\mathbb K,q,\alpha}$,   
$(q,\alpha)\in\I$, 
is linearly disjoint over $\mathbb K$.
\end{enumerate}
\end{thm}

\begin{rem}{\rm 
One of the main motivations for \cite{AF26} was to prove, in connection
with Furstenberg's dimension conjecture (see the discussion in
\cite[Sec.\ 2]{AF26}), that an irrational real number cannot be generated
by finite automata in two multiplicatively independent integer bases
\cite[Thm.\ 2.2]{AF26}. This result  follows immediately from
Theorem~\ref{thm:M-values}. 
Indeed, suppose that a real number $\xi$ admits expansions in the multiplicatively independent integer
bases $b_1$ and $b_2$ whose digit sequences are generated by finite
automata with respective input bases $q_1$ and $q_2$. Then
\[
\xi\in
\M_{\mathbb Q,q_1,1/b_1}
\cap
\M_{\mathbb Q,q_2,1/b_2}.
\]
Since $b_1$ and $b_2$ are multiplicatively independent, the pairs
$(q_1,1/b_1)$ and $(q_2,1/b_2)$ are not equivalent. Hence, by
Theorem~\ref{thm:M-values}, the two algebras above are linearly disjoint
over $\mathbb Q$, and their intersection is therefore equal to
$\mathbb Q$. It follows that $\xi$ is rational.}
\end{rem}

Theorem~\ref{thm:M-values} yields the structural decomposition of the $\Q$-algebra $\M$ announced above.  

\begin{coro}
Let $\mathcal R$ be a set of representatives of 
$
\bigl(\mathbb Z_{\ge2}\times(\Q^\star\cap\mathcal D)\bigr)/\sim$. 
Then the canonical homomorphism of $\Q$-algebras 
\[
\bigotimes_{(q,\alpha)\in\mathcal R}
\M_{\Q,q,\alpha}
\longrightarrow
\M,
\qquad
\bigotimes_{(q,\alpha)\in\mathcal R} x_{q,\alpha}
\longmapsto
\prod_{(q,\alpha)\in\mathcal R} x_{q,\alpha}
\]
is an isomorphism. Here, the tensor product is taken over $\Q$, with
$x_{q,\alpha}=1$ for all but finitely many $(q,\alpha)$ in every
elementary tensor.
\end{coro}

The study of linear and algebraic relations among the elements of $\M$ 
thus reduces to the corresponding problem within each of the smaller 
algebras $\M_{\Q,q,\alpha}$. Within each of these algebras, the lifting theorems
\cite[Th.\ 1.4]{AF17} and \cite[Th.\ 1]{AF24} provide powerful tools
for determining algebraic relations. This does not mean, however, that the remaining problem is
always easy. Indeed, each algebra $\M_{\Q,q,\alpha}$ has infinite
transcendence degree over $\Q$. This follows, for instance, from Mahler's
work \cite{Ma30a}, by considering the values at $\alpha$ of the Fredholm
series\footnote{The terminology was introduced by Kempner \cite{Ke},
following a paper of Mittag-Leffler that referred to work of Fredholm.
The latter work, however, had no connection with the series under
consideration, which would therefore more appropriately be called the
Kempner series (see the discussions in \cite{Ad13,Sh92,Za20}).
Nevertheless, the name \emph{Fredholm series} has become so widespread
that we shall refrain from changing it.}
\begin{equation}\label{eq:Fredholm_def}
\mathfrak f_q(z)=\sum_{n\geq 0}z^{q^n}
\end{equation}
and of its derivatives; see also Section~\ref{sec:F}.

Theorem~\ref{thm:M-values} also implies a linear disjointness result for
$M$-functions, generalizing \cite[Thm.\ 2.4]{AF26} and
\cite[Thm.\ 1.8]{ADHW}\footnote{Theorem~1.8 in \cite{ADHW} is stated
over the base field $\mathbb C$, but the extension from $\Q$ to
$\mathbb C$ follows from a routine argument.}.

\begin{coro}\label{coro:functions}
Let $\mathcal M$ be the $\Q(z)$-algebra generated by all $M$-functions
and, for every $q\in\mathbb N_{\ge2}$, let $\mathcal M_q$ be the
$\Q(z)$-algebra generated by all $M_q$-functions.
Set
$
\mathcal B
=
\{q\ge2:\text{$q$ is not a perfect power}\}$. 
Then the canonical homomorphism of $\Q(z)$-algebras
\[
\bigotimes_{q\in\mathcal B}\mathcal M_q
\longrightarrow
\mathcal M,
\qquad
\bigotimes_{q\in\mathcal B}f_q
\longmapsto
\prod_{q\in\mathcal B}f_q
\]
is an isomorphism.
Here, the tensor product is taken over $\Q(z)$, with $f_q=1$ for all but
finitely many $q$ in every elementary tensor.
\end{coro}

Since the generating series of $q$-automatic sequences are
$M_q$-functions, this corollary generalizes Cobham's celebrated theorem
on automatic sequences (see, for instance, the discussion in
\cite[Sec.\ 2.3]{AF26}). In fact, it yields a substantially stronger
statement: all linear and algebraic relations over $\Q(z)$ among
generating series of automatic sequences are obtained simply by
combining the corresponding relations within each family of
$q$-automatic generating series, where $q$ ranges over the integers
that are not perfect powers.

We also mention another linear disjointness result concerning base field
extensions, recently obtained by the authors in
\cite[Cor.\ 3]{AF24}. If
$\mathbb K\subset\mathbb L\subset\Q$, then
$$
\M_{\mathbb K,q,\alpha}
\otimes_{\mathbb K}\mathbb L
\simeq
\M_{\mathbb L,q,\alpha}.
$$
This descent theorem further reduces the study of algebraic relations
over $\Q$ to the corresponding problem over the smaller field generated
by the coefficients of the underlying power series and the evaluation
point.

Theorem~\ref{thm:M-values} is obtained as a consequence of two recent
results of the authors on multivariate Mahler systems: a purity theorem
and an extension of Nishioka's theorem; see
\cite[Cor.\ 3.5 and Cor.\ 3.9]{AF26}.

\subsection{Application to Fredholm series}\label{sec:F}

As an application of our general results, we give a complete description
of all algebraic relations among the values of the Fredholm series
introduced in \eqref{eq:Fredholm_def} and their derivatives at nonzero algebraic points in the
open unit disk. These series are among the oldest and most classical
examples in Mahler's method. They appeared in Mahler's foundational work
\cite{Ma29,Ma30a,Ma30b} and have since been repeatedly used as basic
examples throughout the subject
\cite{Ku76_a,Ku76_b,Ku77,Mah82,Ni91,Ni_Liv,Ma99,Ma_Liv}. Algebraic
relations among their values have been explicitly studied in
\cite{Ku76_a,LvdP78,Ni94,Ni01,Ni02}. Nevertheless, a complete description involving arbitrary parameters,
algebraic points, and derivatives was still missing. 

There are several evident sources of relations among Fredholm series.
They give rise to the following four rules:
\begin{enumerate}
\item[{\rm (F)}]\label{item:Mahler}
the \emph{functional-equation rule}
$
\mathfrak f_q(\alpha)
=
\mathfrak f_q(\alpha^q)+\alpha,
$
which follows from the $q$-Mahler equation
$
\mathfrak f_q(z)=\mathfrak f_q(z^q)+z;
$

\item[{\rm (R)}]\label{item:roots}
the \emph{root-of-unity rule}
$
\mathfrak f_q(\zeta\alpha)
=
\sum_{i=0}^{\varphi(d)-1}
\zeta^{q^i}
\mathfrak f_{q^{\varphi(d)}}(\alpha^{q^i}),
$
where $d$ is coprime to $q$, $\zeta$ is a primitive $d$th root of
unity, and $\varphi$ denotes Euler's totient function;

\item[{\rm (E)}]\label{item:scale}
the \emph{exponentiation rule}
$
\mathfrak f_q(\alpha)
=
\sum_{i=0}^{n-1}
\mathfrak f_{q^n}(\alpha^{q^i})$, $n\geq 2$, 
which relates values associated with the parameter $q$ to values
associated with the parameter $q^n$;

\item[{\rm (D)}]\label{item:deriv}
the \emph{derivative rule}, obtained by successively differentiating
the functional identities underlying the rules
{\rm (F)}, {\rm (R)}, and {\rm (E)}, and then specializing at algebraic
points.
\end{enumerate}

The exponentiation rule follows by splitting the defining series
according to the residue class of the exponent index modulo $n$. Indeed,
writing $k=nj+i$, with $0\leq i\leq n-1$, we obtain
\begin{equation}\label{eq:E}
\mathfrak f_q(z)
=
\sum_{i=0}^{n-1}\sum_{j\geq0}z^{q^{nj+i}}
=
\sum_{i=0}^{n-1}
\mathfrak f_{q^n}(z^{q^i}).
\end{equation}
Specialization at $z=\alpha$ gives rule {\rm (E)}.
The root-of-unity rule is obtained by applying \eqref{eq:E} with
$z=\zeta\alpha$ and $n=\varphi(d)$, and using
$
\zeta^{q^{i+j\varphi(d)}}=\zeta^{q^i}
$
for  $j\geq0$.

Our main result concerning the Fredholm series and their derivatives
asserts that these evident relations generate all algebraic relations
among the corresponding values. The difficulty is not to produce such
relations, but to prove that no additional relation can arise from the
interaction between distinct Mahler parameters, distinct algebraic points, and different
orders of derivation. We state here an informal version; the precise formulation is given and
proved in Section~\ref{sec:Fredholm}
(see Theorem~\ref{thm:Fredholm}). 

\begin{thm}\label{thm:Fredholm_intro}
Every algebraic relation among values of Fredholm series and their
derivatives at nonzero algebraic points in the open unit disk can be
obtained by successive applications of the rules {\rm (F)}, {\rm (R)},
{\rm (E)}, and {\rm (D)}.
\end{thm}

Thus, Theorem~\ref{thm:Fredholm_intro} provides an explicit description
of the entire ideal of algebraic relations. This completely
resolves a classical problem that has accompanied Mahler's method since
its origins.

\begin{rem}{\rm 
Let us illustrate Theorem~\ref{thm:Fredholm_intro} with an example.
In \cite{Ni02}, Nishioka concludes her paper with the following identity:
\begin{equation}\label{eq:Ni}
\sum_{k\geq 0}\alpha^{3^k}
-\sqrt{-1}\sum_{k\geq 0}
\bigl(\sqrt{-1}\,\alpha\bigr)^{3^k}
-2\sum_{k\geq 0}\alpha^{9^k}=0.
\end{equation}
She observes that this relation illustrates the difficulty of finding
necessary and sufficient conditions for the algebraic independence of
values of Fredholm series. In our notation, relation~\eqref{eq:Ni} reads
\begin{equation}\label{eq:Ni2}
\mathfrak f_3(\alpha)
-\sqrt{-1}\,\mathfrak f_3\bigl(\sqrt{-1}\,\alpha\bigr)
-2\mathfrak f_9(\alpha)=0.
\end{equation}

Applying rule {\rm (E)} with $q=3$ and $n=2$ gives
\begin{equation}\label{eq:Ee}
\mathfrak f_3(\alpha)
=
\mathfrak f_9(\alpha)+\mathfrak f_9(\alpha^3).
\end{equation}
On the other hand, applying rule {\rm (R)} with
$q=3$, $d=4$, and $\zeta=\sqrt{-1}$, we obtain
$$
\mathfrak f_3\bigl(\sqrt{-1}\,\alpha\bigr)
=
\sqrt{-1}\,\mathfrak f_9(\alpha)
-\sqrt{-1}\,\mathfrak f_9(\alpha^3),
$$
since $\varphi(4)=2$,
$\bigl(\sqrt{-1}\bigr)^3=-\sqrt{-1}$, and
$3^{\varphi(4)}=9$. Multiplying the latter identity by $\sqrt{-1}$
and subtracting it from \eqref{eq:Ee}, we obtain
$$
\mathfrak f_3(\alpha)
-\sqrt{-1}\,\mathfrak f_3\bigl(\sqrt{-1}\,\alpha\bigr)
=
2\mathfrak f_9(\alpha),
$$
which is equivalent to \eqref{eq:Ni2}, and hence to \eqref{eq:Ni}. Thus, Nishioka's relation \eqref{eq:Ni} is obtained by combining the elementary rules
{\rm (E)} and {\rm (R)}.}
\end{rem}

\section{Background in Mahler's method}
\label{sec:Mahler}

At the end of the 1920s, Mahler
\cite{Ma29,Ma30a,Ma30b} developed a pioneering method for studying
algebraic relations among values of $M$-functions. Important progress
has recently made it possible to consider simultaneously values of
$M_q$-functions associated with distinct parameters $q$ and evaluated
at several distinct algebraic points; see \cite{AF26,AF26_ad}. We recall
below the results that will be needed in the sequel.

Let $n$ be a positive integer and let
$
T=(t_{i,j})_{1\leq i,j\leq n}
$
be an $n\times n$ matrix with nonnegative integer entries. Given an
$n$-tuple of indeterminates
$
\z=(z_1,\ldots,z_n),
$
we set
$$
T\z
:=
\left(
\prod_{j=1}^n z_j^{t_{1,j}},
\ldots,
\prod_{j=1}^n z_j^{t_{n,j}}
\right),
$$
and we let $T$ act similarly on $(\mathbb C^\star)^n$.

A \emph{linear $T$-Mahler system}, or simply a \emph{Mahler system}, is
a system of functional equations of the form
\begin{equation}\label{eq:T_Mahler_sys}
\left(
\begin{array}{c}
f_1(\z)\\
\vdots\\
f_m(\z)
\end{array}
\right)
=
A(\z)
\left(
\begin{array}{c}
f_1(T\z)\\
\vdots\\
f_m(T\z)
\end{array}
\right),
\end{equation}
where
$
f_1(\z),\ldots,f_m(\z)\in\bQ[[\z]]
$
and
$
A(\z)\in{\rm GL}_m(\bQ(\z)).
$

An algebraic point
$
\balpha\in(\bQ^\star)^n
$
is said to be \emph{regular with respect to} the system
\eqref{eq:T_Mahler_sys} if $A(\z)$ is well defined and invertible at
$T^\ell\balpha$ for every $\ell\in\N$. In several variables, the
application of Mahler's method requires an additional admissibility
condition on the pair $(T,\balpha)$; see \cite[Sec.\ 5]{AF26}.

In the present
paper, however, we shall only need scalar matrices, that is, matrices of
the form
$
T=q{\rm I}_n
$
for some integer $q\geq2$. In this case,
$
T\z=(z_1^q,\ldots,z_n^q)
$
and
$
T(\alpha_1,\ldots,\alpha_n)
=
(\alpha_1^q,\ldots,\alpha_n^q).
$
This restriction allows one to express the admissibility condition on
the pair $(T,\balpha)$ entirely in terms of the point $\balpha$. Recall that $\mathcal D$ denotes the open unit disk. 
More precisely, the pair
$
(T=q{\rm I}_n,\balpha)
$
is admissible in the sense of \cite{AF26} if and only if
$\balpha\in(\Q^\star\cap\mathcal D)^n$ and the coordinates of $\balpha$ are globally
multiplicatively independent, that is,
$$
\alpha_1^{u_1}\cdots\alpha_n^{u_n}=1,
\qquad
(u_1,\ldots,u_n)\in\mathbb Z^n,
$$
implies
$
u_1=\cdots=u_n=0.
$

The following result is a particular case of \cite[Cor.\ 3.5]{AF26}. 
It is a multivariate extension of Nishioka's theorem \cite[Thm.\ 4.2.1]{Ni_Liv}.

\begin{thm}
\label{thm:lifting}
Let $n$ be a positive  integer, $T=q{\rm I}_n$, and  
$
f_1(\z),\ldots,f_m(\z)\in\bQ[[\z]]
$
be related by a $T$-Mahler system. Let
$
\balpha=(\alpha_1,\ldots,\alpha_n)
\in(\bQ^\star\cap\mathcal D)^n.
$
Assume that the coordinates of $\balpha$ are globally multiplicatively
independent and that $\balpha$ is regular with respect to this system.
Then
$$
{\rm tr.deg}_{\bQ}
\bigl(
f_1(\balpha),\ldots,f_m(\balpha)
\bigr)
=
{\rm tr.deg}_{\bQ(\z)}
\bigl(
f_1(\z),\ldots,f_m(\z)
\bigr).
$$
\end{thm}

We shall also need a purity theorem for systems associated with
multiplicatively independent parameters. Let $r$ be a positive integer.
For every $i\in\{1,\ldots,r\}$, consider a $T_i$-Mahler system
\begin{equation}\tag{\theequation.$i$}\label{eq:mahleriT}
\left(
\begin{array}{c}
f_{i,1}(\z_i)\\
\vdots\\
f_{i,m_i}(\z_i)
\end{array}
\right)
=
A_i(\z_i)
\left(
\begin{array}{c}
f_{i,1}(T_i\z_i)\\
\vdots\\
f_{i,m_i}(T_i\z_i)
\end{array}
\right),
\end{equation}
where $n_i$ and $m_i$ are positive integers,
$
\z_i=(z_{i,1},\ldots,z_{i,n_i})
$
is a vector of indeterminates,
$
T_i=q_i{\rm I}_{n_i}
$
with $q_i\geq2$,
$
A_i(\z_i)\in{\rm GL}_{m_i}(\bQ(\z_i)),
$
and
$
f_{i,1}(\z_i),\ldots,f_{i,m_i}(\z_i)
\in\bQ[[\z_i]].
$
Let also
$
\balpha_i\in(\bQ^\star\cap\mathcal D)^{n_i}.
$

The following purity theorem is a particular case of \cite[Cor.\ 3.9]{AF26}. 

\begin{thm} 
\label{thm:purete}
With the above notation, assume that the following conditions hold.
\begin{enumerate}[label=\rm(\roman*)]
\item For every $i\in\{1,\ldots,r\}$, the coordinates of $\balpha_i$
are globally multiplicatively independent, and $\balpha_i$ is regular
with respect to the corresponding $T_i$-Mahler system.

\item The integers $q_1,\ldots,q_r$ are pairwise multiplicatively
independent.
\end{enumerate} 
For every $i\in\{1,\ldots,r\}$, set
$
\mathcal E_i
=
\{
f_{i,1}(\balpha_i),\ldots,f_{i,m_i}(\balpha_i)
\},
$
and set
$
\mathcal E=\bigcup_{i=1}^r\mathcal E_i.
$
Then
$$
{\rm tr.deg}_{\bQ}(\mathcal E)
=
\sum_{i=1}^r
{\rm tr.deg}_{\bQ}(\mathcal E_i).
$$
\end{thm}

\section{Proof of Theorem \ref{thm:M-values}}\label{sec:thm1}

In this section, we prove our main general result concerning values of $M$-functions.

\subsection{Proof of part {\rm (a)} of
Theorem~\ref{thm:M-values}}

We first establish the following elementary stability lemma.

\begin{lem}\label{lem:stab_Mp}
Let $n,q\geq2$ be integers and let $f\in\Q[[z]]$. Then the following
assertions hold.
\begin{enumerate}[label=\rm(\roman*)]
\item\label{item:si_power}
If $f(z^n)$ is an $M_q$-function, then $f(z)$ is an $M_q$-function.

\item\label{item:alors_power}
If $f(z)$ is an $M_q$-function, then $f(z^n)$ is an $M_q$-function.

\item\label{item:decompo}
Suppose that $f(z)$ is an $M_q$-function and write
$
f(z)=\sum_{i=0}^{n-1}z^if_i(z^n),
$
with $f_i(z)\in\Q[[z]]$. Then each $f_i(z)$ is an $M_q$-function.
\end{enumerate}
\end{lem}

\begin{proof}
Let us first prove \ref{item:si_power}. Consider a $q$-Mahler equation
satisfied by $f(z^n)$:
\begin{equation}\label{eq:Mahler_power}
a_0(z)f(z^n)+a_1(z)f(z^{qn})
+\cdots+a_m(z)f(z^{q^mn})=0,
\end{equation}
where $a_0,\ldots,a_m\in\Q[z]$ are not all zero. 
For every
$r\in\{0,\ldots,m\}$, write
$
a_r(z)=\sum_{j=0}^{n-1}z^ja_{r,j}(z^n),
$
with $a_{r,j}(z)\in\Q[z]$. Substituting these decompositions into
\eqref{eq:Mahler_power} and identifying the powers of $z$ modulo $n$,
we obtain, for every $j\in\{0,\ldots,n-1\}$,
$$
a_{0,j}(z^n)f(z^n)
+a_{1,j}(z^n)f(z^{qn})
+\cdots
+a_{m,j}(z^n)f(z^{q^mn})=0.
$$ 
Replacing $z^n$ by $z$, we infer that
$$
a_{0,j}(z)f(z)
+a_{1,j}(z)f(z^q)
+\cdots
+a_{m,j}(z)f(z^{q^m})=0.
$$ 
Since some $a_i$ is nonzero, there exists a pair $(i,j)$ such that $a_{i,j}\neq0$. Hence
$f$ satisfies a nontrivial $q$-Mahler equation and is therefore an
$M_q$-function. This proves \ref{item:si_power}.

Let us now prove \ref{item:alors_power}. Let
$$
a_0(z)f(z)+a_1(z)f(z^q)+\cdots+a_m(z)f(z^{q^m})=0
$$ 
be a nontrivial $q$-Mahler equation satisfied by $f$. Replacing $z$
by $z^n$, we obtain
$$
a_0(z^n)f(z^n)+a_1(z^n)f(z^{qn})
+\cdots+a_m(z^n)f(z^{q^mn})=0.
$$ 
Thus, $f(z^n)$ is an $M_q$-function, proving
\ref{item:alors_power}.

It remains to prove \ref{item:decompo}. Since $f$ is an
$M_q$-function, it occurs as a coordinate of a vector
$F(z)\in\Q[[z]]^m$ satisfying a $q$-Mahler system
$$
F(z)=A(z)F(z^q), \quad A(z)\in{\rm GL}_m(\Q(z)).
$$
Write
$
F(z)=\sum_{i=0}^{n-1}z^iF_i(z^n),
$
with $F_i(z)\in\Q[[z]]^m$. Then
$$
\sum_{i=0}^{n-1}z^iF_i(z^n)
= 
A(z)\sum_{i=0}^{n-1}z^{qi}F_i(z^{qn})=
\sum_{i=0}^{n-1}z^{qi}A(z)F_i(z^{qn}).
$$
Identifying the powers of $z$ modulo $n$, we see that every coordinate
of every vector $F_i(z^n)$ is a linear combination over  $\Q(z)$ of the
coordinates of
$
F_0(z^{qn}),\ldots,F_{n-1}(z^{qn}).
$

Let $G(z)$ be the vector obtained by concatenating
$
F_0(z^n),\ldots,F_{n-1}(z^n).
$
It follows that there exists a matrix
$B(z)\in\mathcal M_{mn}(\Q(z))$ such that
$
G(z)=B(z)G(z^q).
$
By scalar elimination, every coordinate of a vector satisfying such a
linear Mahler system satisfies a scalar $q$-Mahler equation; the matrix
$B(z)$ need not be invertible for this conclusion. Thus, each
$f_i(z^n)$ is an $M_q$-function. Assertion
\ref{item:si_power} then implies that each $f_i(z)$ is an
$M_q$-function.
\end{proof}

\begin{proof}[Proof of part {\rm (a)} of
Theorem~\ref{thm:M-values}]
Suppose first that
$
(q,\alpha)\not\sim(r,\beta).
$
Then either $q$ and $r$ are multiplicatively independent, or $\alpha$
and $\beta$ are multiplicatively independent. Let
$
\xi\in
\M_{\mathbb K,q,\alpha}
\cap
\M_{\mathbb K,r,\beta}.
$
Then
$
\xi=f_1(\alpha)=f_2(\beta),
$
where $f_1\in\mathbb K[[z]]$ is an $M_q$-function and
$f_2\in\mathbb K[[z]]$ is an $M_r$-function. Since these two values are
equal, they are algebraically dependent. It follows from
\cite[Thm.\ 1.1]{AF26} that $\xi\in\mathbb K$. Hence
$
\M_{\mathbb K,q,\alpha}
\cap
\M_{\mathbb K,r,\beta}
=
\mathbb K.
$
Moreover,
$
\M_{\mathbb K,q,\alpha}\neq\mathbb K,
$
since it contains the transcendental number
$\mathfrak f_q(\alpha)$, where $\mathfrak f_q$ is the Fredholm series introduced in \eqref{eq:Fredholm_def}. Therefore,
$
\M_{\mathbb K,q,\alpha}
\neq
\M_{\mathbb K,r,\beta}.
$
This proves one implication.

Suppose now that
$
(q,\alpha)\sim(r,\beta).
$
Then $q$ and $r$ are multiplicatively dependent, and so are $\alpha$
and $\beta$. Since $q,r\geq2$, there exist positive integers $u$ and
$v$ such that
$
q^u=r^v.
$

We first note that, for every positive integer $s$, the classes of
$M_q$-functions and $M_{q^s}$-functions coincide. Indeed, every
$M_q$-function is an $M_{q^s}$-function, by iterating a $q$-Mahler
equation. Conversely, every $q^s$-Mahler equation is also a
$q$-Mahler equation. Consequently,
$
\M_{\mathbb K,q,\alpha}
=
\M_{\mathbb K,q^u,\alpha}
$
and
$
\M_{\mathbb K,r,\beta}
=
\M_{\mathbb K,r^v,\beta}.
$
Replacing $q$ and $r$ by their common power $q^u=r^v$, which we again 
denote by $q$, we are reduced to proving that 
$
\M_{\mathbb K,q,\alpha}
=
\M_{\mathbb K,q,\beta}.
$

Since $\alpha$ and $\beta$ are multiplicatively dependent and satisfy
$0<|\alpha|,|\beta|<1$, there exist positive integers $a$ and $b$ such
that
$
\alpha^a=\beta^b.
$
It is therefore enough to prove that, for every positive integer $t$,
$
\M_{\mathbb K,q,\alpha}
=
\M_{\mathbb K,q,\alpha^t}.
$
Indeed, applying this equality first to $\alpha$ and then to $\beta$, and using
$\alpha^a=\beta^b$, we obtain
$
\M_{\mathbb K,q,\alpha}
=
\M_{\mathbb K,q,\alpha^a}
=
\M_{\mathbb K,q,\beta^b}
=
\M_{\mathbb K,q,\beta}.
$

Let first
$
\xi\in\M_{\mathbb K,q,\alpha^t}.
$
There exists an $M_q$-function $f\in\mathbb K[[z]]$ such that
$
\xi=f(\alpha^t).
$
By Lemma~\ref{lem:stab_Mp}\ref{item:alors_power}, the series
$
g(z):=f(z^t)
$
is an $M_q$-function belonging to $\mathbb K[[z]]$. Since
$
\xi=g(\alpha),
$
we obtain
$
\M_{\mathbb K,q,\alpha^t}
\subseteq
\M_{\mathbb K,q,\alpha}.
$

Conversely, let
$
\xi\in\M_{\mathbb K,q,\alpha}.
$
There exists an $M_q$-function $f\in\mathbb K[[z]]$ such that
$
\xi=f(\alpha).
$
Write
$
f(z)=\sum_{i=0}^{t-1}z^if_i(z^t),
$
where $f_i(z)\in\mathbb K[[z]]$. By
Lemma~\ref{lem:stab_Mp}\ref{item:decompo}, each $f_i$ is an
$M_q$-function. Hence
$
g(z):=\sum_{i=0}^{t-1}\alpha^if_i(z)
$
is also an $M_q$-function in $\mathbb K[[z]]$. By construction,
$
g(\alpha^t)
=
\sum_{i=0}^{t-1}\alpha^if_i(\alpha^t)
=
f(\alpha)
=
\xi.
$
Thus,
$
\M_{\mathbb K,q,\alpha}
\subseteq
\M_{\mathbb K,q,\alpha^t}.
$
This proves
$
\M_{\mathbb K,q,\alpha}
=
\M_{\mathbb K,q,\alpha^t},
$
which completes the proof.
\end{proof}

\subsection{Proof of part {\rm (b)} of
Theorem~\ref{thm:M-values}}

The proof is divided into three steps. First, for a fixed parameter $q$,
we express values of $M_q$-functions at several algebraic points as values
of functions satisfying a multivariate Mahler system at a regular and
admissible point. Next, we use a functional purity theorem to separate the families of
functions attached to pairwise multiplicatively independent points. 
Finally, we combine the multivariate extension of Nishioka's theorem with the purity theorem (Theorem~\ref{thm:purete}) to establish the required algebraic independence statement.

Given
$
\bmu=(\mu_1,\ldots,\mu_n)\in\N^n
$
and
$
\balpha=(\alpha_1,\ldots,\alpha_n)\in\Q^n,
$
we use the standard monomial notation
$
\balpha^\bmu
:=
\alpha_1^{\mu_1}\cdots\alpha_n^{\mu_n}.
$
Similarly, if $\z=(z_1,\ldots,z_n)$ is a vector of indeterminates, we set
$
\z^\bmu:=z_1^{\mu_1}\cdots z_n^{\mu_n}.
$

The following lemma is a straightforward generalization of
\cite[Prop.\ A.7]{AF26_ad}. In that proposition, only one function is
considered at each point $\alpha_i$. The same proof applies to finitely
many functions at each point, by enlarging the corresponding Mahler
systems.

\begin{lem}\label{lem:existenceMpoints}
Let $q\geq2$ be an integer, and let
$\alpha_1,\ldots,\alpha_r\in\Q$ be pairwise multiplicatively independent
algebraic numbers such that
$
0<|\alpha_i|<1
$
for every $i$. For each $i\in\{1,\ldots,r\}$, let
$
f_{i,1}(z),\ldots,f_{i,s_i}(z)\in\Q[[z]]
$
be $M_q$-functions that are well defined at $\alpha_i$.

Then there exist positive integers $t$ and $\ell$, a point
$\bbeta\in\Q^t$, a matrix $T\in\mathcal M_t(\N)$, vectors
$
\bmu_1,\ldots,\bmu_r\in\N^t,
$
and, for each $i\in\{1,\ldots,r\}$, a root of unity $\zeta_i$, an
integer $m_i\geq s_i$, and $M_q$-functions
$
g_{i,1}(z),\ldots,g_{i,m_i}(z)\in\Q[[z]],
$
such that the following assertions hold.
\begin{enumerate}[label=\rm(\alph*)]
\item\label{item:prop-points-a}
For every $i\in\{1,\ldots,r\}$,
$
\alpha_i=\zeta_i\bbeta^{\bmu_i}.
$

\item\label{item:prop-points-b}
For every $i\in\{1,\ldots,r\}$ and 
$j\in\{1,\ldots,s_i\}$,
$
f_{i,j}(\alpha_i)=g_{i,j}(\bbeta^{\bmu_i}).
$

\item\label{item:prop-points-c}
For every $i\in\{1,\ldots,r\}$, the functions
$
g_{i,1}(z),\ldots,g_{i,m_i}(z)
$
are related by a $q^\ell$-Mahler system, and
$\bbeta^{\bmu_i}$ is regular with respect to this system.

\item\label{item:prop-points-d}
The $m_1+\cdots+m_r$ functions
$
g_{i,j}(\z^{\bmu_i})$, 
$
1\leq i\leq r$, $1\leq j\leq m_i$, 
are related by a $T$-Mahler system, where
$\z=(z_1,\ldots,z_t)$ is a vector of indeterminates.

\item\label{item:prop-points-e}
The point $\bbeta$ is regular for this $T$-Mahler system, its
coordinates belong to $\Q^\star\cap \mathcal D$, and they are globally
multiplicatively independent. 

\item\label{item:prop-points-f}
We have
$
T=q^\ell{\rm I}_t.
$

\item\label{item:prop-points-g}
The vectors
$
\bmu_1,\ldots,\bmu_r
$
are pairwise linearly independent over $\Q$.
\end{enumerate}
\end{lem}

The next result is \cite[Lem.\ A.5]{AF26_ad}.

\begin{lem}\label{lem:purity_functions}
Let $r$ and $t$ be positive integers, and let
$
\bmu_1,\ldots,\bmu_r\in\N^t
$
be pairwise linearly independent over $\Q$. For each
$i\in\{1,\ldots,r\}$, let $m_i$ be a positive integer and let
$
f_{i,1}(z),\ldots,f_{i,m_i}(z)\in\Q[[z]].
$
If $\z=(z_1,\ldots,z_t)$ is a vector of indeterminates, then
\begin{multline*}
{\rm tr.deg}_{\Q(\z)}
\bigl(
f_{i,j}(\z^{\bmu_i})
:
1\leq i\leq r,\ 1\leq j\leq m_i
\bigr)
\\
=
\sum_{i=1}^r
{\rm tr.deg}_{\Q(z)}
\bigl(
f_{i,j}(z)
:
1\leq j\leq m_i
\bigr).
\end{multline*}
\end{lem}

We now recall some standard facts concerning linear disjointness and
freeness. The definitions given in the introduction extend the usual
notions for field extensions. Indeed, two $\mathbb K$-algebras
$\mathbb E$ and $\mathbb F$ are linearly disjoint over $\mathbb K$ if
and only if their fields of fractions are linearly disjoint over
$\mathbb K$. Likewise, they are free over $\mathbb K$ if and only if
their fields of fractions are algebraically disjoint over $\mathbb K$,
in the terminology of Bourbaki. We follow here the terminology of Lang;
see \cite[Chap.\ VIII]{Lang} and
\cite[Chap.\ V, Sec.\  1 and~14]{Bourbaki}.

The corresponding statements also hold for arbitrary families of
algebras. Moreover, if the field of fractions of each algebra in the
family is a regular extension of $\mathbb K$, then linear disjointness
and freeness are equivalent; see
\cite[Chap.\ VIII, Thm.\ 4.12]{Lang} and
\cite[Chap.\ V, Sec.\ 14]{Bourbaki}.

Recall that an integral $\mathbb K$-algebra $\mathbb E$ is said to be
regular over $\mathbb K$ if its field of fractions
$\operatorname{Frac}(\mathbb E)$ is a regular extension of $\mathbb K$.
In characteristic zero, this is equivalent to requiring that
$\mathbb K$ be algebraically closed in
$\operatorname{Frac}(\mathbb E)$.

The next result is a rephrasing of \cite[Cor.\ 1.8]{AF17}. 

\begin{lem}
\label{lem:regularity}
Let $\mathbb K\subset\Q$. The algebras
$
\M_{\mathbb K,q,\alpha}$, 
$q\in\mathbb Z_{\geq2}$, $
\alpha\in\mathbb K^\star\cap \mathcal D$,  
are regular over $\mathbb K$.
\end{lem}

We are now ready to prove part {\rm (b)} of
Theorem~\ref{thm:M-values}.

\begin{proof}[Proof of part {\rm (b)} of
Theorem~\ref{thm:M-values}]
By Lemma~\ref{lem:regularity}, the field of fractions of each
algebra $\M_{\mathbb K,q,\alpha}$ is a regular extension of
$\mathbb K$. It follows from
\cite[Chap.\ VIII, Thm.\ 4.12]{Lang} that the algebras
$\M_{\mathbb K,q,\alpha}$, 
$
(q,\alpha)\in\I$, 
are linearly disjoint over $\mathbb K$ if and only if they are free over
$\mathbb K$.

Since freeness is a property of finite subfamilies, we may assume without
loss of generality that $\I$ is finite. Let
$
(q_1,\alpha_1),\ldots,(q_r,\alpha_r)
$
be an enumeration of $\I$.

Using part {\rm (a)} of Theorem~\ref{thm:M-values}, we may replace each
pair $(q_i,\alpha_i)$ by an equivalent pair without changing the
corresponding algebra. We may therefore assume that, for every $i<j$,
either $q_i$ and $q_j$ are multiplicatively independent, or
$
q_i=q_j
$
and $\alpha_i$ and $\alpha_j$ are multiplicatively independent.

For every $i\in\{1,\ldots,r\}$, choose elements
$
\xi_{i,1},\ldots,\xi_{i,s_i}
\in
\M_{\mathbb K,q_i,\alpha_i}
$
that are algebraically independent over $\mathbb K$. Since $\Q/\mathbb K$ is algebraic, they are also algebraically independent
over $\Q$. We prove that the
$
s_1+\cdots+s_r
$
numbers
$
\xi_{i,j}$, 
$
1\leq i\leq r$, $1\leq j\leq s_i$,  
are algebraically independent over $\Q$, and hence over
$\mathbb K$. This will prove that the algebras
$
\M_{\mathbb K,q_i,\alpha_i}$, $1\leq i\leq r$, 
are free over $\mathbb K$.

For every $i\in\{1,\ldots,r\}$ and
$j\in\{1,\ldots,s_i\}$, choose an $M_{q_i}$-function
$
f_{i,j}(z)\in\mathbb K[[z]]
$
such that
$
f_{i,j}(\alpha_i)=\xi_{i,j}.
$ 
Partition $\{1,\ldots,r\}$ into classes
$
\mathcal I_1,\ldots,\mathcal I_t
$
so that $i$ and $j$ belong to the same class $\mathcal I_k$ if and only if
$
q_i=q_j=:\rho_k.
$

For each $k\in\{1,\ldots,t\}$, the points
$
\alpha_i$, 
$i\in\mathcal I_k$, 
are pairwise multiplicatively independent. We may therefore apply
Lemma~\ref{lem:existenceMpoints} with $q=\rho_k$ to these points
and to the corresponding functions
$
f_{i,1}(z),\ldots,f_{i,s_i}(z).
$
Let
$
\bbeta_k$, 
$(\bmu_i)_{i\in\mathcal I_k}$, 
$\ell_k$, 
$T_k$, 
$\z_k$, 
$(g_{i,j}(z))$ ($i\in\mathcal I_k$, $1\leq j\leq m_i$),
and $B_k(\z_k)$ 
denote the point, vectors, iteration level, transformation matrix, vector
of indeterminates, family of functions, and matrix of the resulting
$T_k$-Mahler system.

Lemma~\ref{lem:existenceMpoints} ensures that
$\bbeta_k$ is regular for the $T_k$-Mahler system
defined by $B_k(\z_k)$, that its coordinates belong to
$\Q^\star\cap\mathcal D$, and that they are globally
multiplicatively independent. Moreover,
$
T_k=\rho_k^{\ell_k}{\rm I}.
$
Since $\rho_1,\ldots,\rho_t$ are pairwise multiplicatively independent,
so are
$
\rho_1^{\ell_1},\ldots,\rho_t^{\ell_t}.
$
We may therefore apply Theorem~\ref{thm:purete} to these $t$ Mahler
systems.

For every $k\in\{1,\ldots,t\}$, set
$
\mathcal F_k
:=
\{
g_{i,j}(\bbeta_k^{\bmu_i})
:
i\in\mathcal I_k,\ 1\leq j\leq m_i
\}.
$
Theorem~\ref{thm:purete} gives
$$
{\rm tr.deg}_{\Q}
\left(
\bigcup_{k=1}^t\mathcal F_k
\right)
=
\sum_{k=1}^t
{\rm tr.deg}_{\Q}(\mathcal F_k).
$$ 

Fix $k\in\{1,\ldots,t\}$. Applying
Theorem~\ref{thm:lifting} to the $T_k$-Mahler system associated with
$B_k(\z_k)$, we obtain
$$
{\rm tr.deg}_{\Q}(\mathcal F_k)
=
{\rm tr.deg}_{\Q(\z_k)}
\bigl(
g_{i,j}(\z_k^{\bmu_i})
:
i\in\mathcal I_k,\ 1\leq j\leq m_i
\bigr).
$$
Since the vectors
$
\bmu_i$, 
$i\in\mathcal I_k$, 
are pairwise linearly independent over $\Q$,
Lemma~\ref{lem:purity_functions} yields
\begin{multline*}
{\rm tr.deg}_{\Q(\z_k)}
\bigl(
g_{i,j}(\z_k^{\bmu_i})
:
i\in\mathcal I_k,\ 1\leq j\leq m_i
\bigr)
\\
=
\sum_{i\in\mathcal I_k}
{\rm tr.deg}_{\Q(z)}
\bigl(
g_{i,j}(z)
:
1\leq j\leq m_i
\bigr).
\end{multline*}

For every $i\in\mathcal I_k$, Condition
\ref{item:prop-points-c} of Lemma~\ref{lem:existenceMpoints}
allows us to apply Theorem~\ref{thm:lifting} to the Mahler system
satisfied by
$
g_{i,1}(z),\ldots,g_{i,m_i}(z)
$
at the regular point $\bbeta_k^{\bmu_i}$. Thus,
$$
{\rm tr.deg}_{\Q(z)}
\bigl(
g_{i,j}(z)
:
1\leq j\leq m_i
\bigr)
=
{\rm tr.deg}_{\Q}(\mathcal F_{k,i}),
$$ 
where
$
\mathcal F_{k,i}
:=
\{
g_{i,j}(\bbeta_k^{\bmu_i})
:
1\leq j\leq m_i
\}.
$ 
Let $\mathcal F = \bigcup_{k=1}^t \mathcal F_k$. Combining the preceding four identities, we obtain
\begin{equation}\label{eq:transdsum}
{\rm tr.deg}_{\Q}(\mathcal F)
=\sum_{k=1}^t
\sum_{i\in\mathcal I_k}
{\rm tr.deg}_{\Q}(\mathcal F_{k,i}).
\end{equation}

For every $i\in\mathcal I_k$, set
$
\mathcal E_{k,i}
:=
\{
f_{i,1}(\alpha_i),\ldots,f_{i,s_i}(\alpha_i)
\}
$
 and 
 $$
 \mathcal E
 :=
 \{
 f_{i,j}(\alpha_i)
 :
 1\leq i\leq r,\ 1\leq j\leq s_i
 \}
 $$
By Condition~\ref{item:prop-points-b} of
Lemma~\ref{lem:existenceMpoints}, we have
$
\mathcal E_{k,i}\subseteq\mathcal F_{k,i}.
$ 
Moreover,
$$
\mathcal F
=\bigcup_{k=1}^t
\bigcup_{i\in\mathcal I_k}\mathcal F_{k,i}
,\quad\text{ and }\quad
\mathcal E
=\bigcup_{k=1}^t
\bigcup_{i\in\mathcal I_k}\mathcal E_{k,i}.
$$
It follows from \cite[Lem.\ 10.3]{AF26} and \eqref{eq:transdsum} that
\begin{eqnarray*}
	{\rm tr.deg}_{\Q}(\mathcal E)
	&=&\sum_{k=1}^t
	\sum_{i\in\mathcal I_k}
	{\rm tr.deg}_{\Q}(\mathcal E_{k,i})\\
	&=&\sum_{i=1}^r
	{\rm tr.deg}_{\Q}
	\bigl(
	f_{i,1}(\alpha_i),\ldots,f_{i,s_i}(\alpha_i)
	\bigr).
\end{eqnarray*} 
For every $i$, the numbers
$
f_{i,1}(\alpha_i),\ldots,f_{i,s_i}(\alpha_i)
$
are algebraically independent over $\Q$, and hence 
$
{\rm tr.deg}_{\Q}
\bigl(
f_{i,1}(\alpha_i),\ldots,f_{i,s_i}(\alpha_i)
\bigr)
=
s_i.
$
Thus,
$$
{\rm tr.deg}_{\Q}(\mathcal E)
=
s_1+\cdots+s_r.
$$ 
It follows that all the numbers
$
\xi_{i,j}$, 
$1\leq i\leq r$, $1\leq j\leq s_i$, 
are algebraically independent over $\Q$, and therefore over
$\mathbb K$. Hence the algebras
$
\M_{\mathbb K,q_i,\alpha_i}$, 
$1\leq i\leq r$,
are free over $\mathbb K$, and thus linearly disjoint over $\mathbb K$. 
This ends the proof. 
\end{proof}

\section{Proof of Corollary~\ref{coro:functions}}

Corollary~\ref{coro:functions} is equivalent to saying that the $\Q(z)$-algebras $\mathcal M_q$, where $q \in \mathcal B$, are linearly disjoint over $\Q(z)$. By \cite[Lem.\ 3.2]{AF17}, each algebra $\mathcal M_q$ is a regular extension of $\Q(z)$. It follows that the algebras $\mathcal M_q$, $q \in \mathcal B$, are linearly disjoint over $\Q(z)$ if and only if they are free over $\Q(z)$. 

Let $q_1,\ldots,q_r \in \mathcal B$ be pairwise distinct and, for every $i \in\{1,\ldots,r\}$, let $f_{i,1},\ldots,f_{i,s_i} \in \mathcal M_{q_i}$ be algebraically independent over $\Q(z)$. We only have to prove that the $s_1+\cdots + s_r$ series $f_{i,j}$, $1\leq i \leq r$, $1\leq j \leq s_i$, are algebraically independent over $\Q(z)$.

Let $R < 1$ be a positive real number. By Case (M$_q$) of \cite[Prop.\ 14]{AF24}, for each $i$, there exists a finite set $\mathcal S_i$ such that the following holds: for any $\alpha \in \bQ^\star \setminus \mathcal S_i$ with $\vert \alpha \vert < R$, the numbers $f_{i,1}(\alpha),\ldots,f_{i,s_i}(\alpha)$ are algebraically independent.

Since each set $\mathcal S_i$ is finite, there exists $\alpha \in \bQ^\star \setminus \bigcup_i \mathcal S_i$ with $\vert \alpha \vert < R$. For such an $\alpha$ and for every $i$, the $s_i$ numbers $f_{i,1}(\alpha),\ldots,f_{i,s_i}(\alpha)$ are algebraically independent. These numbers belong to the algebra $\M_{\Q,q_i,\alpha}$. Since the $q_1,\ldots,q_r \in \mathcal B$ are pairwise distinct, the pairs $(q_i,\alpha)$, $i \in\{1,\ldots,r\}$ are pairwise inequivalent with respect to $\sim$. It follows from part {\rm (b)} of Theorem~\ref{thm:M-values}, applied with $\mathbb K = \bQ$ and $\I=\{(q_1,\alpha),\ldots,(q_r,\alpha)\}$, that the algebras $\M_{\Q,q_1,\alpha},\ldots,\M_{\Q,q_r,\alpha}$ are linearly disjoint over $\Q$. In particular, the $s_1+\cdots+s_r$ numbers
$$
f_{i,j}(\alpha), \quad 1\leq i \leq r,\, 1\leq j \leq s_i
$$
are algebraically independent over $\Q$. A fortiori, the $s_1+\cdots+s_r$ series $f_{i,j}$, $1\leq i \leq r$, $1\leq j \leq s_i$, are algebraically independent over $\Q(z)$.
\qed
\section{Algebraic relations among values of Fredholm series}
\label{sec:Fredholm}

The goal of this section is to give a precise formulation of
Theorem~\ref{thm:Fredholm_intro} and to prove it. To simplify the
exposition, we use the logarithmic derivation, or Euler operator,
$
\theta:=z\frac{\partial}{\partial z},
$
instead of the usual derivation. For every nonnegative integer $\ell$,
we set
$$
\mathfrak f_q^{(\ell)}(z)
:=
\theta^\ell\bigl(\mathfrak f_q\bigr)(z)
=
\sum_{n=0}^{\infty}q^{\ell n}z^{q^n}.
$$

We introduce the $\Q$-algebra generated by all values of Fredholm series
and their successive Euler derivatives at nonzero algebraic points in
the open unit disk:
$$
\mathcal F
:=
\Q\left[
\mathfrak f_q^{(\ell)}(\alpha)
:
q\in\mathbb Z_{\geq2},\
\ell\in\mathbb Z_{\geq0},\
\alpha\in\Q^\star\cap\mathcal D
\right]
\subseteq\mathbb C.
$$
We also introduce the polynomial ring over $\Q$ in the corresponding
countable family of indeterminates:
$$
\mathcal A
:=
\Q\left[
X_{q,\ell,\alpha}
:
q\in\mathbb Z_{\geq2},\
\ell\in\mathbb Z_{\geq0},\
\alpha\in\Q^\star\cap\mathcal D
\right].
$$

We define the evaluation map
\begin{eqnarray*}
{\rm ev}:\mathcal A&\longrightarrow&\mathcal F,\\
X_{q,\ell,\alpha}&\longmapsto&
\mathfrak f_q^{(\ell)}(\alpha).
\end{eqnarray*}
By definition, ${\rm ev}$ is a surjective homomorphism of
$\Q$-algebras.

Let
$
\X=(X_{q,\ell,\alpha})_{q,\ell,\alpha}
$
denote the corresponding family of indeterminates, so that
$
\mathcal A=\Q[\X].
$
We consider the following polynomials:
\begin{equation}\label{eq:poly_fredholm}
\begin{array}{rcl}
F_{q,\ell,\alpha}(\X)
&=&
X_{q,\ell,\alpha}
-q^\ell X_{q,\ell,\alpha^q}
-\alpha,
\\[1ex]
R_{\zeta,q,\ell,\alpha}(\X)
&=&
X_{q,\ell,\zeta\alpha}
-\displaystyle\sum_{i=0}^{\varphi(d)-1}
\zeta^{q^i}q^{i\ell}
X_{q^{\varphi(d)},\ell,\alpha^{q^i}},
\\[2ex]
E_{n,q,\ell,\alpha}(\X)
&=&
X_{q,\ell,\alpha}
-\displaystyle\sum_{i=0}^{n-1}
q^{i\ell}X_{q^n,\ell,\alpha^{q^i}}.
\end{array}
\end{equation}
Here, $q,n,d\geq2$, $\ell\in\mathbb Z_{\geq0}$,
$\alpha\in\Q^\star\cap\mathcal D$, the integers $d$ and $q$ are
coprime, and $\zeta$ is a primitive $d$th root of unity.

When $\ell=0$, these polynomials correspond respectively to the
relations arising from the rules {\rm (F)}, {\rm (R)}, and {\rm (E)}
introduced in Section~\ref{sec:F}. For $\ell\geq1$, they correspond to
the relations obtained by applying $\theta^\ell$ to the underlying
functional identities and then specializing at $\alpha$, in accordance
with rule {\rm (D)}.

We let $\mathcal I$ denote the ideal of $\mathcal A$ generated by all
the polynomials in \eqref{eq:poly_fredholm}. Thus, $\mathcal I$ is the
ideal generated by the relations obtained from the rules
{\rm (F)}, {\rm (R)}, {\rm (E)}, and {\rm (D)}.

Theorem~\ref{thm:Fredholm_intro} can now be stated precisely as
follows.

\begin{thm}\label{thm:Fredholm}
We have
$
\ker({\rm ev})=\mathcal I.
$
Equivalently, the evaluation map induces an isomorphism of
$\Q$-algebras
$$
\overline{{\rm ev}}:
\mathcal A/\mathcal I
\xrightarrow{\;\sim\;}
\mathcal F.
$$
\end{thm}

The proof of Theorem~\ref{thm:Fredholm} is divided into four steps.
First, we show that
$
\mathcal I\subseteq\ker({\rm ev}).
$
To establish the reverse inclusion, we first apply classical results from Mahler's method in one variable to exhibit families of values that are algebraically independent. Using the polynomials in \eqref{eq:poly_fredholm}, we deduce that any relation involving only pairs $(q,\alpha)$ that are pairwise equivalent belongs to $\I$. Then, we leverage the linear disjointness of the algebras $\M_{\Q,q,\alpha}$ to show that the latter relations generate the entire ideal $\ker({\rm ev})$. We conclude that $\ker({\rm ev})= \I$.

\subsection{Elements of $\I$ are relations}

We begin by showing that the polynomials in
\eqref{eq:poly_fredholm} belong to the kernel of ${\rm ev}$.

\begin{lem}\label{lem:Poly_in_Ker}
We have
$
\I\subseteq\ker({\rm ev}).
$
\end{lem}

We recall the following two elementary identities:
\[
\theta^\ell(f(\xi z))=\theta^\ell(f)(\xi z)
\quad\text{and}\quad
\theta^\ell(f(z^r))=r^\ell\theta^\ell(f)(z^r),
\]
where $f$ is a power series, $r$ is a positive integer, and $\xi$ is a complex number.

\begin{proof}
We first consider the polynomials $F_{q,\ell,\alpha}$. Applying
$\theta^\ell$ to the $q$-Mahler equation
$
\mathfrak f_q(z)=\mathfrak f_q(z^q)+z
$
and then specializing at $\alpha$, we obtain
$$
{\rm ev}\bigl(F_{q,\ell,\alpha}\bigr)
=
\mathfrak f_q^{(\ell)}(\alpha)
-q^\ell\mathfrak f_q^{(\ell)}(\alpha^q)
-\alpha
=
0.
$$
Thus,
$
F_{q,\ell,\alpha}\in\ker({\rm ev}).
$

We next consider the polynomials
$
R_{\zeta,q,\ell,\alpha},
$
where $\zeta$ is a primitive $d$th root of unity with $d$ and $q$ relatively prime.
Applying $\theta^\ell$ to the identity
$$
\mathfrak f_q(\zeta z)
=
\sum_{i=0}^{\varphi(d)-1}
\zeta^{q^i}
\mathfrak f_{q^{\varphi(d)}}(z^{q^i})
$$
and then specializing at $\alpha$, we obtain
$$
{\rm ev}\bigl(R_{\zeta,q,\ell,\alpha}\bigr)
=
\mathfrak f_q^{(\ell)}(\zeta\alpha)
-
\sum_{i=0}^{\varphi(d)-1}
\zeta^{q^i}q^{i\ell}
\mathfrak f_{q^{\varphi(d)}}^{(\ell)}(\alpha^{q^i})
=
0.
$$
Therefore,
$
R_{\zeta,q,\ell,\alpha}\in\ker({\rm ev}).
$

Finally, we consider the polynomials $E_{n,q,\ell,\alpha}$. Applying
$\theta^\ell$ to the identity
$
\mathfrak f_q(z)
=
\sum_{i=0}^{n-1}
\mathfrak f_{q^n}(z^{q^i})
$
and then specializing at $\alpha$, we obtain
$$
{\rm ev}\bigl(E_{n,q,\ell,\alpha}\bigr)
=
\mathfrak f_q^{(\ell)}(\alpha)
-
\sum_{i=0}^{n-1}
q^{i\ell}
\mathfrak f_{q^n}^{(\ell)}(\alpha^{q^i})
=
0.
$$
Hence,
$
E_{n,q,\ell,\alpha}\in\ker({\rm ev}).
$

Since $\ker({\rm ev})$ is an ideal of $\mathcal A$, it contains the
ideal generated by all the polynomials in \eqref{eq:poly_fredholm}.
Consequently,
$
\I\subseteq\ker({\rm ev}).
$
\end{proof}
\subsection{Algebraic independence of subfamilies of values}
Fix an integer $q \geq 2$ and an algebraic number $\alpha \in \Q^\star \cap \mathcal D$.
Define an equivalence relation over $\Z_{>0}$, setting
\begin{equation}\label{eq:sim_q}
b\sim_q c
\quad\Longleftrightarrow\quad
\frac{b}{c}=q^a
\quad\text{for some }a\in\mathbb Z.
\end{equation}
Let $(c_\nu)_{\nu\geq 0} \subset \Z_{> 0}$ denote an enumeration of a set of representatives of the equivalence classes with respect to $\sim_q$.
\begin{lem}
	\label{lem:alg_indep_remaining_values}
	The numbers	$\mathfrak f_{q}^{(\ell)}(\alpha^{c_\nu})$, $\ell \geq 0$, $\nu\geq 0$, are algebraically independent. 
\end{lem}
\begin{proof}
	It is sufficient to prove that, for any integers $s,t\geq 0$, the $st$ numbers 
	$$
	\mathfrak f_{q}^{(\ell)}(\alpha^{c_\nu}), \quad 0 \leq \ell < s,\, 0 \leq \nu < t,
	$$
	are algebraically independent. Set
	$$
	g_{\ell,\nu}(z)
	:=
	\mathfrak f_q^{(\ell)}(z^{c_\nu}), \quad 0 \leq \ell < s,\, 0 \leq \nu < t.
	$$
	They satisfy the following $q$-Mahler equation:
	$$
	g_{\ell,\nu}(z)
	=
	q^\ell g_{\ell,\nu}(z^q)+z^{c_\nu}.
	$$
	Hence the function $1$ together with the $st$ functions
	$
	g_{\ell,\nu}(z)$, $0 \leq \ell < s$, $0 \leq \nu < t$,
	satisfies a $q$-Mahler system of size $st+1$ whose matrix has polynomial entries and nonzero constant
	determinant.
	In
	particular, $\alpha$ is regular with respect to this system. 
	
	We claim that the functions
	$
	g_{\ell,\nu}(z)$, 
	$0 \leq \ell < s$, $0 \leq \nu < t$,
	are algebraically independent over $\Q(z)$. Otherwise,
	\cite[Satz~1]{Ma30b}, or equivalently
	\cite[Cor.\ 9]{Ku77}, would imply the existence of algebraic
	numbers $\mu_{\ell,\nu}$, not all zero, such that
	$$
	\sum_{\ell=0}^{s-1}
	\sum_{\nu=0}^{t-1}
	\mu_{\ell,\nu}g_{\ell,\nu}(z)
	\in\Q(z).
	$$
	This power series is equal to 
	$$
	\sum_{\nu=0}^{t-1}
	\sum_{n\geq0}
	\left(
	\sum_{\ell=0}^{s-1}
	\mu_{\ell,\nu}q^{\ell n}
	\right)
	z^{c_\nu q^n}.
	$$
	Since the integers $c_\nu$ belong to distinct equivalence classes with respect to
	$\sim_q$, the sets
	$
	\{c_\nu q^n:n\geq0\}
	$, $0 \leq \nu < t$,
	are pairwise disjoint. Their union has arbitrarily long gaps. A
	rational power series supported on such a set must therefore be a
	polynomial. It follows that, for every $\nu$,
	$$
	\sum_{\ell=0}^{s-1}
	\mu_{\ell,\nu}q^{\ell n}=0
	$$
	for all sufficiently large $n$. A Vandermonde argument gives
	$
	\mu_{\ell,\nu}=0
	$
	for  every $\ell\in \{0,\ldots,s-1\}$ and every $\nu \in\{0,\ldots,t-1\}$, a contradiction.
	
	Then, Nishioka's theorem, \textit{i.e.}\ Theorem~\ref{thm:lifting} with $n=1$, applied at $\alpha$ implies that
		\begin{align*}
				\label{eq:tr.deg_Fredholm}
					{\rm tr.deg}_\Q((\mathfrak f_q^{(\ell)}(\alpha^{c_\nu}))_{0\leq \ell <s,\, 0 \leq \nu < t}) & =
				{\rm tr.deg}_\Q(1,(g_{\ell,\nu}(\alpha))_{0\leq \ell <s,\, 0 \leq \nu < t}) \\ &= 	{\rm tr.deg}_{\Q(z)}(1,(g_{\ell,\nu}(z))_{0\leq \ell <s,\, 0 \leq \nu < t})
				\\ & = st\,.
			\end{align*}
	Therefore, the $st$ numbers
$
\mathfrak f_q^{(\ell)}(\alpha^{c_\nu})$,  $0 \leq \ell < s$, $0 \leq \nu < t$,
are algebraically independent over $\Q$. Note that the full strength of Nishioka's theorem is not needed here.
Alternatively, one could use the main result of Mahler's paper~\cite{Ma30b},
with $n=1$.
\end{proof}

\subsection{Relations within each equivalence class}

Recall that the equivalence relation $\sim$ on pairs
$
(q,\alpha)\in\mathbb Z_{\geq2}\times(\Q^\star\cap\mathcal D)
$
was defined by
$$
(q,\alpha)\sim(p,\beta)
\quad\Longleftrightarrow\quad
q\sim p
\ \text{and}\
\alpha\sim\beta,
$$
where, on the right-hand side, $\sim$ denotes multiplicative
dependence of non-zero complex numbers.

Let
$
(\mathcal C_i)_{i\in\N}
$
denote an enumeration of the equivalence classes for this relation. For every
$i\in\N$, let $\X_i$ be the subfamily of indeterminates
$X_{q,\ell,\alpha}$ such that
$
(q,\alpha)\in\mathcal C_i.
$
\begin{lem}\label{lem:Fr_basis_lin}
	For every $i \in \N$, we have $\ker({\rm ev}) \cap \Q[\X_i] \subset \I$.
\end{lem}

\begin{proof}
	Let
	$
	P_0\in  \ker({\rm ev}) \cap \Q[\X_i].
	$
	There exist integers $q_1,\ldots,q_s\geq2$, nonnegative integers
	$\ell_1,\ldots,\ell_s$, and algebraic numbers
	$
	\alpha_1,\ldots,\alpha_s\in\Q^\star\cap\mathcal D
	$
	such that
	$$
	P_0\in
	\Q[
	X_{q_1,\ell_1,\alpha_1},
	\ldots,
	X_{q_s,\ell_s,\alpha_s}
	].
	$$
	Since all the pairs $(q_k,\alpha_k)$ belong to the same equivalence
	class $\mathcal C_i$, the integers $q_1,\ldots,q_s$ are pairwise
	multiplicatively dependent, and the numbers
	$\alpha_1,\ldots,\alpha_s$ are pairwise multiplicatively dependent.
	
	There exist an integer $q\geq2$ and positive integers
	$u_1,\ldots,u_s$ such that
	$$
	q_k=q^{u_k},
	\qquad
	1\leq k\leq s.
	$$
	There also exist an algebraic number
	$
	\alpha\in\Q^\star\cap\mathcal D,
	$
	roots of unity $\zeta_1,\ldots,\zeta_s$, and positive integers
	$v_1,\ldots,v_s$ such that
	$$
	\alpha_k=\zeta_k\alpha^{v_k},
	\qquad
	1\leq k\leq s.
	$$
	
	We first remove from each $\zeta_k$ the part of its order supported on
	the prime divisors of $q$. 
	The relation
	$$
	F_{q_k,\ell_k,\beta}
	=
	X_{q_k,\ell_k,\beta}
	-q_k^{\ell_k}X_{q_k,\ell_k,\beta^{q_k}}
	-\beta
	$$
	shows that, modulo $\I$, the indeterminate
	$
	X_{q_k,\ell_k,\beta}
	$
	is equal to
	$$
	q_k^{\ell_k}X_{q_k,\ell_k,\beta^{q_k}}+\beta.$$ 
	Applying this operation sufficiently many times, we obtain
	$$
	P_0=P_1+Q_0,
	\qquad
	Q_0\in \I ,
	$$
	where $P_1$ involves only indeterminates of the form
	$
	X_{q_k,\ell_k,\eta_k\alpha^{a_k}}$, 
	with $a_k\geq1$ and $\eta_k$ a root of unity whose order $d_k$ is
	coprime with $q$, and hence with $q_k$.
	
	The relation
	$$
	R_{\eta_k,q_k,\ell_k,\alpha^{a_k}}
	=
	X_{q_k,\ell_k,\eta_k\alpha^{a_k}}
	-
	\sum_{h=0}^{\varphi(d_k)-1}
	\eta_k^{q_k^h}q_k^{h\ell_k}
	X_{q_k^{\varphi(d_k)},\ell_k,\alpha^{a_kq_k^h}}
	$$
	then allows us to eliminate all roots of unity. Thus,
	$$
	P_1=P_2+Q_1,
	\qquad
	Q_1\in \I,
	$$
	where $P_2$ involves only indeterminates of the form
	$
	X_{r_k,\ell_k,\alpha^{a_k}}$, 
	and the integers $r_k$ are powers of $q$.
	
	Choose a common power $r$ of the integers $r_k$. Thus, for every $k$,
	there exists a positive integer $n_k$ such that
	$
	r_k^{n_k}=r.
	$
	Using the relation
	$$
	E_{n_k,r_k,\ell_k,\alpha^{a_k}}
	=
	X_{r_k,\ell_k,\alpha^{a_k}}
	-
	\sum_{h=0}^{n_k-1}
	r_k^{h\ell_k}
	X_{r,\ell_k,\alpha^{a_kr_k^h}},
	$$
	we obtain
	$$
	P_2=P_3+Q_2,
	\qquad
	Q_2\in \I,
	$$
	where $P_3$ involves only
	indeterminates of the form
	$
	X_{r,\ell,\alpha^b}$, $\ell\geq 0$, $b\geq1$.
	
	Consider an enumeration $(c_\nu)_{\nu \geq 0}$ of a set of representatives for the equivalence relation $\sim_r$ introduced in \eqref{eq:sim_q} such that $r \nmid c_\nu$, for every $\nu$. Then, each element in the equivalence class of $c_\nu$ is of the form
	$
	c_\nu r^a
	$
	with $a\geq0$. Set
	$
	\beta_\nu:=\alpha^{c_\nu}.
	$
	Every indeterminate appearing in $P_3$ is of the form $X_{r,\ell,\alpha^{c_\nu r^a}}=X_{r,\ell,\beta_\nu^{r^a}}$, for some $\nu,a\geq 0$. 
	
	Repeated use of the relations
	$$
	F_{r,\ell,\beta_\nu^{r^{a-1}}}
	=
	X_{r,\ell,\beta_\nu^{r^{a-1}}}
	-r^\ell X_{r,\ell,\beta_\nu^{r^{a}}}
	-\beta_\nu^{r^{a-1}},\quad a \geq1,
	$$
	allows us to replace every indeterminate
	$
	X_{r,\ell,\beta_\nu^{r^a}}
	$
	by an expression involving
	$
	X_{r,\ell,\beta_\nu}
	$
	and an algebraic constant. Thus, we may replace every indeterminate
	$
	X_{r,\ell,\alpha^b}
	$ 
	by an expression involving
	$
	X_{r,\ell,\alpha^{c_\nu}}
	$
	for some $\nu \geq 0$. We therefore obtain
	$$
	P_3=P_4+Q_3,
	\qquad
	Q_3\in \I,
	$$
	where $P_4$ involves only indeterminates of the form $X_{r,\ell,\alpha^{c_\nu}}$, $\ell \geq 0$, $\nu\geq 0$. 
	
	We have the decomposition
	$$
	P_0 = Q_0 + Q_1 + Q_2 + Q_3 + P_4.
	$$
	By Lemma~\ref{lem:Poly_in_Ker}, $\I \subset \ker({\rm ev})$, so $Q_0$, $Q_1$, $Q_2$, and $Q_3$ all belong to $\ker({\rm ev})$. Since $P_0 \in \ker({\rm ev})$ as well, we have $P_4 \in \ker({\rm ev})$. The polynomial $P_4$ involves only indeterminates of the form $X_{r,\ell,\alpha^{c_\nu}}$ with $\ell \geq 0$ and $\nu \geq 0$. Lemma~\ref{lem:alg_indep_remaining_values} implies that the numbers $\mathfrak f_{r}^{(\ell)}(\alpha^{c_\nu})$ are algebraically independent. Thus, $P_4=0$ and
	$$
	P_0 = Q_0 + Q_1 + Q_2 + Q_3 \in \I.
	$$
	This proves the inclusion $\ker({\rm ev}) \cap \Q[\X_i]\subset \I$.
	
\end{proof}
\subsection{Linear disjointness}
We continue with the notation $\mathcal C_i$ and $\X_i$ from the previous subsection.

\begin{lem}\label{lem:Fredholm_Lin_Dis}
The ideal $\ker({\rm ev})$ is generated by the elements of the sets
$
\ker({\rm ev})\cap \Q[\X_i]$, 
$
i\in\N$. 
\end{lem}

\begin{proof}
For every $i\in\N$, choose a pair
$
(q,\alpha)\in\mathcal C_i
$
and set
$
\M_i:=\M_{\Q,q,\alpha}.
$
By Theorem~\ref{thm:M-values}\ref{item:equal}, this algebra does not
depend on the choice of $(q,\alpha)$ in $\mathcal C_i$.

If
$
(q,\alpha)\in\mathcal C_i
$
and
$
\ell\in\mathbb Z_{\geq0},
$
then
$$
{\rm ev}(X_{q,\ell,\alpha})
=
\mathfrak f_q^{(\ell)}(\alpha)
\in\M_i,
$$
since $\mathfrak f_q^{(\ell)}$ is an $M_q$-function. Hence the image of
the restriction
$$
{\rm ev}_i:
\Q[\X_i]
\longrightarrow
\mathbb C
$$
is a subalgebra of $\M_i$. 
Let $\mathcal K$ be a finite subset of $\N$, and set
$$
\mathcal A_{\mathcal K}
:=
\Q\left[
\bigcup_{i\in\mathcal K}\X_i
\right].
$$
Since the families $\X_i$, $i\in\mathcal K$, are pairwise disjoint,
there is a natural isomorphism
$$
\mathcal A_{\mathcal K}
\simeq
\bigotimes_{i\in\mathcal K}\Q[\X_i],
$$
where the tensor products are taken over $\Q$. 
For every $i\in\mathcal K$, set
$$
\mathcal J_i
:=
\ker({\rm ev})\cap\Q[\X_i].
$$
We obtain a natural surjective homomorphism
$$
\mathcal A_{\mathcal K}
\longrightarrow
\bigotimes_{i\in\mathcal K}
\Q[\X_i]/\mathcal J_i,
$$
whose kernel is the ideal of $\mathcal A_{\mathcal K}$ generated by
the elements of $\mathcal J_i$, $i\in\mathcal K$. 
Since the classes $\mathcal C_i$, $i\in\mathcal K$, are distinct,
Theorem~\ref{thm:M-values} implies that the algebras
$
\M_i,
$
$i\in\mathcal K$, are linearly disjoint over $\Q$. Consequently, their
subalgebras
$
{\rm ev}_i\bigl(\Q[\X_i]\bigr),
$
$i\in\mathcal K$, are also linearly disjoint over $\Q$. It follows that
the multiplication map
$$
\bigotimes_{i\in\mathcal K}
\Q[\X_i]/\mathcal J_i
\longrightarrow
\mathbb C
$$
induced by evaluation is injective. Therefore,
$
\ker({\rm ev})\cap\mathcal A_{\mathcal K}
$
is precisely the ideal of $\mathcal A_{\mathcal K}$ generated by the
elements of $\mathcal J_i$, $i\in\mathcal K$.

Finally, every polynomial in $\mathcal A$ involves indeterminates from
only finitely many families $\X_i$. Taking $\mathcal K$ to be the
corresponding finite set of indices proves that $\ker({\rm ev})$ is
generated by the elements of the sets
$
\ker({\rm ev})\cap\Q[\X_i]$, for
$i\in\N$.
\end{proof}

\subsection{End of the proof of Theorem~\ref{thm:Fredholm}}

By Lemma~\ref{lem:Poly_in_Ker}, we already know that
$
\I\subseteq\ker({\rm ev}).
$
It remains to prove the reverse inclusion.

Lemma~\ref{lem:Fredholm_Lin_Dis} shows that $\ker({\rm ev})$ is generated by the elements of the sets $\ker({\rm ev})\cap \Q[\X_i]$, for $i \in \N$. Furthermore, by Lemma~\ref{lem:Fr_basis_lin}, $\ker({\rm ev}) \cap\Q[\X_i] \subset \I$ for every $i \in \N$. Therefore,
$\ker({\rm ev}) \subseteq \I,$
which completes the proof of Theorem~\ref{thm:Fredholm}.
\qed

\subsection*{Declaration on the use of AI}
All new mathematical ideas and results presented in this article
are entirely due to the authors. AI tools were used solely
to check spelling and identify typographical errors in
mathematical expressions.


\end{document}